\documentclass[11pt]{article}
\PassOptionsToPackage{obeyspaces}{url}
\usepackage{amsfonts, amsmath, amssymb}

\usepackage[hidelinks]{hyperref}
\usepackage{bookmark}
\bookmarksetup{open,numbered,depth=5} % paragraphs go into outlines

\usepackage{amsthm}

\usepackage{breakurl}

\usepackage{etoolbox} 
\patchcmd{\thebibliography}{\leftmargin\labelwidth}{\leftmargin\labelwidth\addtolength\itemsep{-0.3\baselineskip}}{}{}

\usepackage{mathtools}

\usepackage{enumitem}

\usepackage[nameinlink]{cleveref}
\crefname{enumi}{}{}

\newtheorem{theorem}{Theorem}
\newtheorem{lemma}[theorem]{Lemma}
\newtheorem{corollary}[theorem]{Corollary}
\newtheorem{proposition}[theorem]{Proposition}
\newtheorem{claim}[theorem]{Claim}
\newtheorem{observation}[theorem]{Observation}
\crefname{conjecture}{Conjecture}{Conjectures}
\crefname{observation}{Observation}{Observations}
\crefname{proposition}{Proposition}{Propositions}
\crefname{claim}{Claim}{Claims}
\newcommand*\samethanks[1][\value{footnote}]{\footnotemark[#1]} % From https://tex.stackexchange.com/questions/4170/multiple-thanks-that-refer-to-same-text
\author{Boris Bukh\thanks{Department of Mathematical Sciences, Carnegie Mellon University, Pittsburgh, PA 15213, USA\@. Supported in part by U.S.\ taxpayers through NSF grants DMS-2154063 and DMS-2452120.}\,\,\textsuperscript{,}\thanks{Supported in part by a Simons Foundation Fellowship.} \and Aleksandre Saatashvili\samethanks[1]}

\title{Most frequent subsequences in a word}

\date{}

\newcommand*{\eqdef}{\stackrel{\mbox{\normalfont\tiny def}}{=}}   % definition by equality                                      
\newcommand*{\N}{\mathbb{N}}                                     % Natural numbers
\newcommand*{\E}{\mathop{\mathbb{E}}}                            % Expectation
\newcommand*{\Z}{\mathbb{Z}}                                     % Integers
\DeclarePairedDelimiter\abs{\lvert}{\rvert}                     % Absolute values, cardinality
\DeclareMathOperator{\Shape}{Sh}                               % Shape of a subsequence

\newcommand*{\LCS}{\text{LCS}}

\usepackage{stmaryrd}
\DeclarePairedDelimiter\ind{\llbracket}{\rrbracket}                    % letter at particular index

\begin{document}
\maketitle
\begin{abstract}
    We prove that every $n$-letter word over $k$-letter alphabet contains some word as a subsequence
    in at least $k^{n/4k(1+o(1))}$ many ways, and that this is sharp as $k\to\infty$. For fixed $k$,
    we show that the analogous number deviates from $\mu_k^n$, for some constant $\mu_k$, by
    a factor of at most~$n$.
\end{abstract}
\section*{Introduction}
A word $v$ is a \emph{subsequence} of a word $w$ if $v$ can be obtained from $w$ by deleting some symbols in $w$.  There might be several ways to obtain $v$ as a subsequence of $w$. For example, $\textsc{abra}$ is a subsequence of the word $\textsc{abracadabra}$ in 9 ways,
which include $\underline{\textsc{abr}}\textsc{ac}\underline{\textsc{a}}\textsc{dabra}$ and $\underline{\textsc{ab}}\textsc{racadab}\underline{\textsc{ra}}$ among others.

Let $M(v,w)$ be the number of occurrences of $v$ as a subsequence in the word~$w$, and let
\[
  M(w)\eqdef \max_v M(v,w)
\]
be the frequency of the most common subsequence.

Consider all $n$-letter words over the alphabet $[k]\eqdef \{1,2,\dotsc,k\}$, for $k\geq 2$, and define
\[
  M_k(n)\eqdef \min_{w\in [k]^n} M(w).
\]
Fang \cite{fang2025maximalnumbersubwordoccurrences} noted that the limit
$\mu_k\eqdef \lim_{n\to\infty} M_k(n)^{1/n}$ exists and proved that $\mu_2\leq\sqrt{1+\sqrt{2}}$ and 
$\mu_k\geq 1+1/k$. 

Our first result is about the rate of convergence of $M_k(n)^{1/n}$ to $\mu_k$.
\begin{theorem}\label[theorem]{Rate thm}
  Suppose that $k\geq 2$ and $n\geq 3$ are integers. Then
  \[
    M_k(n)^{1/n}\leq \mu_k\leq \bigl(n M_k(n)\bigr)^{1/n}.
  \]
\end{theorem}
Here, the lower bound on $\mu_k$ is from \cite{fang2025maximalnumbersubwordoccurrences}, and the upper bound is new.
The theorem gives a way to compute, in principle, $\mu_k$ to arbitrary precision.

Our principal result is that $(\mu_k)^k$ grows like $k^{1/4}$.
\begin{theorem}\label[theorem]{Very Main}
    There exists an absolute constant $C>0$, such that, for $n\geq Ck\log k$, $$\frac
    {1}{C \log^2 k}\cdot k^{1/4}\leq M_k(n)^{k/n}\leq Ck^{1/4}.$$
\end{theorem}

\section{Basic definitions}
For words $w_1$ and $w_2$, the notation $w_1w_2$ denotes their concatenation. For a word $w$ and an integer~$x$, the notation $w^x$ denotes the concatenation of $x$ copies of $w$.

If $w \in [k]^n$ is a word, we write $|w| \eqdef n$ for its length.

For a word $w$, we write $w\ind{i}$ to denote its $i$-th symbol. For an interval $[a,b]$, we use the notation $v\ind{a,b}\eqdef v\ind{a} v\ind{a+1}\dotsb v\ind{b}$ for the subword of $v$ composed of symbols in positions indexed by the interval.

A \emph{map} from a word $v\in [k]^m$ to another word $w\in [k]^n$ is a strictly increasing function $f\colon [m]\to[n]$, such that $v \ind{i}=w \ind{f(i)}$. Let $\mathcal{M}(v,w)$  be the set of all maps from $v$ to $w$. Note that
\[
  M(v,w)=\abs{\mathcal{M}(v,w)}.
\]

\section{Rate of convergence for \texorpdfstring{$\mu_k$}{μ\textunderscore k}}\label{rate sec}
As we mentioned, Fang \cite{fang2025maximalnumbersubwordoccurrences} noted that the inequality $M(n)^{1/n}\leq \mu_k$ follows from the supermultiplicativity
property $M(ww')\geq M(w)M(w')$. The same inequality implies existence of the limit $\mu_k$ via Fekete's lemma. To bound the rate of convergence,
we prove that the function $M_k(n)$ is approximately submultiplicative.

To illustrate the idea, we begin by proving a slightly simpler result first.
\begin{proposition}\label[proposition]{submult lemma} % The optional parameter of \label is needed to counter the buggy version of cleveref that arXiv uses
For all positive integers $m,n$, we have
\[
  M_k(mn)\leq \binom{mn+m-1}{m-1} M_k(n)^m.
\]
\end{proposition}
This implies $\mu_k\leq \bigl(en M_k(n)\bigr)^{1/n}$ by raising both sides to the power $\frac{1}{mn}$, and letting~$m\to\infty$.
\begin{proof}[Proof of \Cref{submult lemma}]
Let $w\in [k]^n$ be any word of length $n$, and consider the $m$-th power of~$w$. Pick
a word $v$  such that $M(v,w^m)\geq M_k(mn)$.
Such a word exists by the definition of the function $M_k(mn)$.
Let $L\eqdef \abs{v}$ be its length.

For $f\in M(v,w^m)$, let $I_j(f)\eqdef \{i\in [L] : (j-1)n+1\leq f(i)\leq jn\}$ for $j\in [m]$.
Evidently, $I_j(f)$ is the largest interval whose $f$-image is contained in the $j$'th copy of $w$ inside~$w^m$.
Write $I(f)\eqdef\nobreak \bigl(I_1(f),\dotsc,I_m(f)\bigr)$, and define $\mathcal{M}_I\eqdef \{f\in \mathcal{M}(v,w^m) : I(f)=I\}$.

Fix a collection of intervals $I=(I_1,\dotsc,I_m)$, and let 
\[
  \mathcal{N}_j(I)\eqdef \{ f|_{I_j} : f\in \mathcal{M}_{I} \},
\]
and note the inclusion $\mathcal{N}_j(I)\subset \mathcal{M}(v\ind{I_j},w)$ for all $j\in [m]$. These inclusions
and the fact that the function $f\mapsto(f|_{I_1},f|_{I_2},\dotsc,f|_{I_m})$ is an injection from
$\mathcal{M}_I$ to $\mathcal{N}_1(I)\times \mathcal{N}_2(I)\times\dotsb\times \mathcal{N}_m(I)$, together imply
that
\begin{equation}\label{mi bound}
  \abs{\mathcal{M}_I}\leq \prod_{j=1}^m \abs{\mathcal{N}_j(I)}\leq \prod_{j=1}^m M(v\ind{I_j},w)\leq \prod_{\ell=0}^n M(w\mid \ell)^{S(I;\ell)},
\end{equation}
where $S(I\mid \ell)$ is the number of intervals of length exactly $\ell$ among $I_1,\dotsc,I_m$, and
\begin{equation}\label{m ell}
  M(w\mid \ell)\eqdef\max_{v\in [k]^\ell} M(v,w).
\end{equation}

In particular, this implies
\begin{equation}\label{basic mi}
  \abs{\mathcal{M}_I}\leq M(w)^{\sum_{\ell}S(I;\ell)}= M(w)^m.
\end{equation}

Since there are $\binom{L+m-1}{m-1}\geq \binom{mn+m-1}{m-1}$ partitions
of the interval $[L]$ into $m$ subintervals, applying inequality~\eqref{basic mi} to each such partition
yields $M_k(mn)\leq \mathcal{M}(v,w^m)\leq \binom{mn+m-1}{m-1} M(w)^m$.
\end{proof}

\paragraph{Proof of the upper bound in \texorpdfstring{\Cref{Rate thm}}{Theorem 1}.}
The bound in the following theorem is stronger than that in \Cref{submult lemma} when $m$ grows faster than~$n$.
\begin{theorem}\label{thm technical rate}
  Suppose that $w\in [k]^n$ is an $n$-letter word, and let $M(w\mid \ell)$ be as in \eqref{m ell}.
  Then
  \[
    M_k(mn)\leq \binom{m+n}{n}\left(\sum_{\ell=0}^n M(w\mid \ell)\right)^{m}.
  \]
\end{theorem}
From $M(w,0)=M(w,n)=1$ and $M_k(n)\geq M(w\mid 1)\geq \lceil n/2\rceil\geq 2$ it follows that $\sum_{\ell=0}^n M(w\mid \ell)\leq (n-1)M_k(n)+2\leq nM_k(n)$ follows. So
\Cref{thm technical rate} result implies the upper bound in \Cref{Rate thm} by taking $m\to\infty$.
\begin{proof} We proceed as in the proof of \Cref{submult lemma}, except we bound the right
  side of \eqref{mi bound} more carefully than the naive \eqref{basic mi}.

  Let
  \[
    \mathcal{R}\eqdef \{(r_0,r_1,\dotsc,r_n)\in \Z^{n+1} : r_i\geq 0,\ \sum_{\ell=0}^n r_\ell=m\}.
  \]
  Recall that the multinomial coefficient are $\binom{m}{r_0,\dotsc,r_n}\eqdef \frac{m!}{r_0!\dotsb r_n!}$,
  and that every $r\in \mathcal{R}$ satisfies
  \begin{equation}\label{multinom ineq}
    1=\left(\frac{r_0}{m}+\dotsb+\frac{r_n}{m}\right)^{m}\geq \left(\frac{r_0}{m}\right)^{r_0}\cdot \ldots\cdot \left(\frac{r_n}{m}\right)^{r_n} \binom{m}{r_0,\dotsc,r_n}.
  \end{equation}
  % TODO: Make sure that the page breaks are OK despite the huge align*
  \allowdisplaybreaks
  
  We compute
  \begin{alignat*}{2}
    M_k(mn)&\leq \mathcal{M}(v,w^m)=\sum_I \abs{\mathcal{M}_I}\\
           &\leq \sum_I \prod_{\ell=0}^n M(w\mid \ell)^{S(I;\ell)}&\qquad&\text{using \eqref{mi bound}}\\
           &=\sum_{r\in \mathcal{R}} \binom{m}{r_0,r_1,\dotsc,r_n}\prod_{\ell=0}^n M(w\mid \ell)^{r_{\ell}}\\
           &\leq \sum_{r\in \mathcal{R}} \prod_{\ell=0}^n \left(\frac{m}{r_\ell}\right)^{r_\ell} M(w\mid \ell)^{r_\ell}&\qquad&\text{using \eqref{multinom ineq}}\\
           &= \sum_{r\in \mathcal{R}} \left(\prod_{\ell=0}^n \left(\frac{m M(w\mid\ell)}{r_\ell}\right)^{r_\ell/m}\right)^m\\
           &\leq \sum_{r\in \mathcal{R}} \left(\sum_{\ell=0}^n M(w\mid\ell)\right)^m&\qquad&\text{using the AM-GM inequality}\\
           &= \binom{m+n}{n} \left(\sum_{\ell=0}^n M(w\mid\ell)\right)^m.&&\tag*{\qedhere}
% The \tag* trick is from https://tex.stackexchange.com/questions/256528/how-to-set-qedhere-properly-within-alignat-environment
   \end{alignat*}
\end{proof}
\paragraph{A computation.} Fang \cite{fang2025maximalnumbersubwordoccurrences} computed the values of $M_2(n)$ for $n\leq 40$. In particular,
he found that $M_2(40)=5500610$, which implies $1.474\leq \mu_2\leq 1.617$ via \Cref{Rate thm}.
If instead one applies \Cref{thm technical rate} directly to the same $40$-symbol word from \cite{fang2025maximalnumbersubwordoccurrences}
that shows the bound $M_2(40)\leq 5500610$, one obtains $\mu_2\leq 1.566$. These bounds do not improve on $\tfrac{3}{2}\leq \mu_2\leq \sqrt{1+\sqrt{2}}\approx 1.554$
from \cite{fang2025maximalnumbersubwordoccurrences}.

\section{Proof of the lower bound in \texorpdfstring{\Cref{Very Main}}{Theorem 2}}
Since the function $M_k(n)$ is monotone as a function of $k$, it suffices to prove
the lower bound $M_k(n)^{k/n}\geq k^{1/4}/C\log^2 k$ only for $k$ on some growing sequence in which the ratios
between the consecutive terms are bounded.
With hindsight, we choose $k\eqdef 2^r-(2^r\bmod r^2)$ for $r\in\mathbb{N}$. We assume that $r$ is large enough.
These choices ensure that $2^{r-1}<k\leq 2^r$ and that $k/r^2$ is an integer.

\begin{observation}\label[observation]{concat}
For any words $w_1,w_2$ we have  $M(w_1w_2)\geq M(w_1)M(w_2)$.
\end{observation}
\begin{proof}
This follows from $M(s_1s_2,w_1w_2)\geq M(s_1,w_1)M(s_2,w_2)$.
\end{proof}
\begin{observation}\label[observation]{identical}
For any word $w$, we have $M(ww)\geq |w|$.
\end{observation}
\begin{proof}
Indeed, $M(w,ww)\geq |w|+1$.
\end{proof}
A \emph{permutation} of a set $\Sigma\subset [k]$ is a word of length $\abs{\Sigma}$ in which each symbol of $\Sigma$ occurs exactly once.
For words $w_1,w_2$ we denote by $\LCS(w_1,w_2)$ the length of the longest common subsequence between $w_1$ and $w_2$.
The following is a version of a lemma due to Beame and Huynh-Ngoc \cite[Lemma 4]{bh}
(with an identical proof, which we include for completeness).
\begin{lemma}\label[lemma]{lcslemma}
  Let $\pi_1,\pi_2,\pi_3$ be permutations over $\Sigma$. Then
  $\LCS(\pi_1,\pi_2)\LCS(\pi_1,\pi_3)\LCS(\pi_2,\pi_3)\geq \abs{\Sigma}$.
\end{lemma}
\begin{proof} For a symbol $\sigma\in \Sigma$, let $L_{i,j}(\sigma)$ be the length
  of the longest common subsequence between $\pi_i$ and $\pi_j$ that ends in~$\sigma$.
  Set $L(\sigma)\eqdef \bigl(L_{1,2}(\sigma),L_{1,3}(\sigma),L_{2,3}(\sigma)\bigr)$.

  The function $L\colon \Sigma\to \N^3$
  is injective. Indeed, suppose that $\sigma,\sigma'\in \Sigma$ are distinct symbols.
  There are two permutations in which they appear in the same order, say $\sigma$ precedes $\sigma'$ in both $\pi_i$
  and $\pi_j$. Then $L_{i,j}(\sigma)<L_{i,j}(\sigma')$, which implies $L(\sigma)\neq L(\sigma')$.

  Since $L$ takes values in $[\LCS(\pi_1,\pi_2)]\times [\LCS(\pi_1,\pi_3)]\times [\LCS(\pi_2,\pi_3)]$, the lemma follows from injectivity of $L$.
\end{proof}

\begin{lemma}\label[lemma]{lcs-kst}
  Suppose that each of $\pi_1,\dotsc,\pi_{2r^2+r}\in [k]^{k/r^2}$ is a permutation on some (possibly different) $(k/r^2)$-element subset of $[k]$.
  Then, for some triple $i<j<\ell$, we have
\begin{equation}\label{eq-prodlcs}
  \LCS(\pi_i,\pi_j)\LCS(\pi_i,\pi_{\ell})\LCS(\pi_j,\pi_{\ell})\geq k/r^{8}.
\end{equation}
\end{lemma}
\begin{proof}
  Define an auxiliary bipartite graph $G$ on parts $[2r^2+r]$ and $[k]$, where we connect $s\in \Sigma$ with
  $i\in [2r^2+r]$ if $\pi_i$ contains the letter~$s$. Since the degree of each vertex in $[2r^2+r]$ is $k/r^2$,
  the average degree of vertices in $[k]$ is $2+1/r$. Pick $\{i,j,\ell\}$ uniformly at random from
  $\binom{[2r^2+r]}{3}$. Let $N$ be the common neighborhood of these three vertices.
  By the usual convexity argument applied to the function
  that is equal to $\binom{x}{3}$ for $x\geq 2$, we deduce that $\E[\abs{N}]$
  is minimized when $k/r$ vertices in $[k]$ have degree $3$, and the rest have degree~$2$.
  So, $\E[\abs{N}]\geq k/r\binom{2r^2+r}{3}$.

  Therefore, there is a choice of $i,j,\ell$ such that $\pi_i,\pi_j,\pi_{\ell}$ have at least
  $k/r\binom{2r^2+r}{3}\geq k/r^8$ symbols in common. An application of \Cref{lcslemma} to this triple completes
  the proof.
\end{proof}

We are now ready to prove our main lemma.
\begin{lemma}\label[lemma]{main}
For all $w\in [k]^{2k+5k/r}$ we have $M(w)\geq k^{1/2}/2r^4$.
\end{lemma}
\begin{proof}
  Partition $w$ into $2r^2+5r$ consecutive subwords $w_1,\dotsc,w_{2r^2+5r}$ of length $k/r^2$ each.
  Let $R\eqdef \{i\in [2r^2+5r] : w_i\text{ is a permutation}\}$. Because $M(w_i)\geq 2$ for $i\notin R$,
  it follows that $M(w)\geq 2^{\abs{[2r^2+5r]\setminus R}}$. So, we may assume that $\abs{R}\geq 2r^2+4r$, for otherwise $M(w)\geq 2^r\geq k$.
  
  Apply \Cref{lcs-kst} to the permutations $w_i$ indexed by $R$ to find a triple $i<j<\ell$
  satisfying \eqref{eq-prodlcs}. Remove the elements $i,j,\ell$ from $R$, and repeat.
  This way we may find a collection $\mathcal{T}$ of disjoint triples of size $\abs{\mathcal{T}}=r$.

  Order the triples in $\mathcal{T}$ according to the middle element of the triple.
  Let $\{i_m<j_m<\ell_m\}$ be the $m$'th triple in this ordering, so that $j_1<\dotsb<j_r$.
  Define
    \begin{align*}
      \alpha_m\eqdef \LCS(w_{i_m},w_{j_m}),\\
      \beta_m\eqdef \LCS(w_{i_m},w_{\ell_m}),\\
      \gamma_m\eqdef \LCS(w_{j_m},w_{\ell_m}).
    \end{align*}    
    Obviously, $M(w)\geq \gamma_1$, $M(w)\geq \alpha_r$, and $M(w)\geq \beta_m$ for all $m\in [r]$.
    In addition, \Cref{concat} implies $M(w)\geq \alpha_m\gamma_{m+1}$ for all $m\in [r-1]$.
    Multiplying these inequalities yields
    \begin{align}
      \label{M-product}
      M(w)^{2r+1}&\geq \prod_{m\in [r]} \alpha_m\beta_m\gamma_m\stackrel{\eqref{eq-prodlcs}}{\geq} (k/r^8)^r,\\
      \intertext{and so}
      \notag M(w)&\geq k^{-1/4r}\cdot k^{1/2}/r^{4}\geq  k^{1/2}/2r^{4}.\qedhere
    \end{align}
\end{proof}
\begin{corollary}\label[corollary]{main cor.}
For a word $w\in [k]^n$ we have $M(w)\geq (ck/\log^8 k)^{n/4k}$, for an absolute constant $c>0$.
\end{corollary}
\begin{proof}
  Let $r$ be as in \Cref{main}. Cut $w$ into $n/(2k+5k/r)$ consecutive subwords of length $2k+5k/r$ each.
  The result follows from \Cref{main} applied to these subwords and \Cref{concat}.
\end{proof}

\section{Proof of the upper bound in \texorpdfstring{\Cref{Very Main}}{Theorem 2}}
Since the function $M_k(n)$ is monotone as a function of $k$, it suffices to prove
the upper bound $M_k(n)^{k/n}\leq (Ck)^{1/4}$ only for $k$ of the form
$k=t^8$. Similarly, since $M_k(n)$ is monotone as a function of $n$, we may
restrict to the values of $n$ that are divisible by $8k$. 

\subsection{Definition and properties of permutations \texorpdfstring{$\pi_i$}{π\textunderscore i}}
\paragraph{Definition.}
Instead of the alphabet $[k]$, we shall use $[t]^8$ as the $k$-letter alphabet in our construction. 

To each \( u \in \{+,-\}^8 \), associate a permutation \( \pi(u) \) on \( [t]^8 \) as follows: The symbol $a\in [t]^8$ precedes $b\in [t]^8$ in $\pi(u)$ if $ua$ is lexicographically smaller than $ub$, where $ux\eqdef (u_1x_1,\dotsc,u_8x_8)\in \{-t,\dotsc,t\}^8$.

%Given \( a, b \in [t]^8 \), let \( i \) be the smallest index such that \( a_i \neq b_i \). Then the symbol \( a \) precedes \( b \) in \( \pi(u) \) if and only if either
%\begin{itemize}
%    \item \( a_i < b_i \) and \( u_i = + \), or
%    \item \( a_i > b_i \) and \( u_i = - \).
%\end{itemize}

The key to our construction are the following eight vectors in $\{-,+\}^8$:
\begin{align*}
u^{(1)}&=(+,+,+,+,+,+,+,+),\\
u^{(2)}&=(-,-,-,+,-,+,-,-),\\
u^{(3)}&=(+,+,+,-,-,-,-,+),\\
u^{(4)}&=(-,+,-,-,+,+,+,-),\\
u^{(5)}&=(+,-,-,+,-,-,+,+),\\
u^{(6)}&=(-,+,+,+,+,-,-,-),\\
u^{(7)}&=(+,-,-,-,+,+,-,+),\\
u^{(8)}&=(-,-,+,-,-,-,+,-).
\end{align*}
Extend them periodically to the sequence $u^{(1)},u^{(2)},\dots$ by setting
$u^{(i)}\eqdef u^{(i\bmod 8)}$, and define
\[
  \pi_i\eqdef \pi(u^{(i)}).
\]
We shall show that the word
\begin{equation}\label{def of w}
  w\eqdef \pi_1\pi_2\dotsc \pi_{n/k}
\end{equation}
obtained by concatenating these permutations satisfies $M(w)^{k/n}\leq (Ck)^{1/4}$. To prove this, we first need to establish several basic properties of the permutations $\pi_i$. 

\paragraph{Properties.}
For $u\in \{-,+\}^8$ and $j\in [8]$, let $u_{\leq j}\eqdef (u_1,u_2,\dotsc,u_j)$ be the projection
onto the first $j$ coordinates. Similarly, define $u_{\geq j}\eqdef (u_j,u_{j+1},\dotsc,u_8)$. We use the same notation $p_{\leq j}$ and $p_{\geq j}$ also for $p\in [t]^8$.

We say that the set of sign vectors $U\subset \{-,+\}^8$ \emph{agrees} at the coordinate $j$ if $u_j$ is the same for all $u\in U$. More generally, we say
that the vectors in $U$ agree on $J\subset [8]$ if they agree at each $j\in J$.

The following properties can be checked by brute force\footnote{The code to verify these
is available at \url{https://www.borisbukh.org/code/maxwords25.html}}.

\begin{proposition}\label[proposition]{Properties of Permutations 1}
The sign vectors $u^{(1)},u^{(2)},\dotsc$ defined above satisfy the following properties:
\begin{enumerate}[label=\alph*), ref=Proposition \thetheorem(\alph*)]
    \item $u^{(i)},u^{(i+1)}$ agree on at most two coordinates.
    \item $u^{(i)},u^{(j)}$ agree on at most four coordinates if $u^{(i)}\neq u^{(j)}$.
    \item $u^{(i)},u^{(i+1)},u^{(i+2)}$ do not agree on any coordinate. 
    \item $u^{(i)},u^{(i+1)},u^{(j)}$ agree on at most one coordinate if $u^{(j)}\notin \{u^{(i)},u^{(i+1)}\}$.
    \item $u^{(i)},u^{(j)},u^{(\ell)}$ agree on at most two coordinates if $u^{(i)}, u^{(j)}, u^{(\ell)}$ are distinct.
    \item \label{prop:u prefix 6} $u^{(i)}_{\geq 7}\neq u^{(i+j)}_{\geq 7}$ for $j\in [3].$
    \item $u^{(i)}_{\geq 6}\neq u^{(i+j)}_{\geq 6}$ for $j\in [7]$.
    \item Vectors $u^{(i)}_{\geq 4},u^{(i+j)}_{\geq 4}$ agree on at most three coordinates for $j\in [7]$.
\end{enumerate}
\end{proposition}
Above we defined $\pi(u)$ for $u\in \{-,+\}^8$. The definition extends
naturally to vectors $u\in \{-,+\}^r$ of any length~$r$; in this case, $\pi(u)$ is a permutation of $[t]^r$. The following lemma is a slight generalization of arguments in \cite{bh,bukh2015twinswordslongcommon}.
\begin{lemma}\label[lemma]{Lemm: intermediate properties}
    Suppose that $U \subset \{+,-\}^r$, and that $J$ is the set
    of positions in which the vectors in $U$ agree. 
    Then $\LCS(\{\pi(u)\}_{u\in U})=t^{\abs{J}}$.
\end{lemma}
\begin{proof}
    Any two symbols of $[t]^r$ that differ only in coordinates in $J$ are ordered in the same manner by all $\pi(u)$, for $u\in U$. So, the lower bound is clear.

    In the opposite direction, suppose that $x$ is a subsequence of $\pi(u)$, for all $u\in U$. If $|x|\geq t^{|J|}+1$, then $x$ contains two distinct symbols, say $x\ind{a}$ and $x\ind{b}$, such that $x\ind{a}_j=x\ind{b}_j$ for all $j\in J$. Let $\ell$ be any coordinate
    on which $x\ind{a}$ and $x\ind{b}$ disagree. Since $\ell\notin J$, there are $u,u'\in U$
    such that $u_\ell\neq u_{\ell}'$. The permutations $\pi(u)$ and $\pi(u')$
    order $x\ind{a}$ and $x\ind{b}$ differently.
    This contradicts $x$ being a subsequence of both $\pi(u)$ and $\pi(u')$. Hence $|x|\leq t^{|J|}$.
\end{proof}

We use \Cref{Lemm: intermediate properties} to translate the properties of sign vectors from \Cref{Properties of Permutations 1} to those of the permutations $\pi_1,\pi_2,\dotsc$.
We write $v\prec w$ to mean that the word $v$ is a subsequence of $w$; similarly, $v \prec w,w'$ means that~$v$ is a common subsequence of words $w,w'$.
\begin{proposition}\label[proposition]{Properties of Permutations 2}
The sequence of permutations $\pi_1,\pi_2,\dotsc$ satisfies the following properties:
\begin{enumerate}[label=\alph*), ref=Proposition \thetheorem(\alph*)]
\item \label{prop: pi adj} $\LCS(\pi_i,\pi_{i+1})\leq t^2$.
\item \label{prop: pi nonadj} $\LCS(\pi_i,\pi_j)\leq t^4$, for $\pi_i\neq \pi_j$.
\item \label{prop: pi three adj} $\LCS(\pi_i,\pi_{i+1},\pi_{i+2})=1$.
\item \label{prop: pi three partadj} $\LCS(\pi_i,\pi_{i+1},\pi_j)\leq t$ if $\pi_j\notin \{\pi_i,\pi_{i+1}\}$.
\item \label{prop:pi distinct} $\LCS(\pi_i,\pi_j,\pi_{\ell})\leq t^2$ if $\pi_i, \pi_j, \pi_\ell$ are distinct.
\item \label{prop:prefix 6} If $j\in [3]$ and $x\prec \pi_i,\pi_{i+j}$ is such that $x\ind{\ell}_{\leq 6}$ is same for all $\ell$, then $\abs{x}\leq t$.
\item \label{prop:prefix 5} If $j\in [7]$ and $x\prec \pi_i,\pi_{i+j}$ is such that $x\ind{\ell}_{\leq 5}$ is same for all $\ell$, then $|x|\leq t^2.$
\item \label{prop:prefix 3} If $j\in [7]$ and $x\prec \pi_i,\pi_{i+j}$ is such that $x\ind{\ell}_{\leq 3}$ is same for all $\ell$, then $|x|\leq t^3.$
\end{enumerate}
\end{proposition}
\begin{proof}
    The properties (a) through (e) follow immediately from the corresponding properties of the $u$-vectors in \Cref{Properties of Permutations 1} and \Cref{Lemm: intermediate properties}.

    To prove (f), define the word $x'$ over the alphabet $[t]^2$ by
    $x'\ind{i}\eqdef x\ind{i}_{\geq 7}$.  By the assumption in~(f), the relations of $x\ind{\ell_1}$ and $x\ind{\ell_2}$ in the lexicographic ordering is
    the same as that between $x'\ind{\ell_1}$ and $x'\ind{\ell_2}$, for any $\ell_1,\ell_2$; and same holds even
    after multiplying $x$ by $u$, and $x'$ by $u_{\geq 7}$, for any $u\in \{-,+\}^8$.
    This    
    implies that $x'$ is a common subsequence of both $\pi(u^{(i)}_{\geq 7})$ and $\pi(u^{(i+j)}_{\geq 7})$. Since $x$ and $x'$ are of the same length,
    an application of \Cref{Lemm: intermediate properties} to the set $U=\{u^{(i)}_{\geq 7},u^{(i+j)}_{\geq 7}\}$ in $\{-,+\}^2$  then concludes the proof via \Cref{prop:u prefix 6}.

%    To prove (f), define the word $x'$ over the alphabet $[t]^2$ by
%    $x'[i]\eqdef x[i]_{\geq 7}$, and note that the assumption in (f) implies that $x'$ is a common subsequence of both $\pi(u^{(i)}_{\geq 7})$ and $\pi(u^{(i+j)}_{\geq 7})$. An application of \Cref{Lemm: intermediate properties} to the set $U=\{u^{(i)}_{\geq 7},u^{(i+j)}_{\geq 7}\}$ in $\{-,+\}^2$  then concludes the proof via \Cref{prop:u prefix 6}.
        
    Similarly, to deduce properties (g) and (h) from the corresponding parts of \Cref{Properties of Permutations 1} we apply \Cref{Lemm: intermediate properties} to $\{u^{(i)}_{\geq 6},u^{(i+j)}_{\geq 6}\}$,  
    and to $\{u^{(i)}_{\geq 4},u^{(i+j)}_{\geq 4}\}$, respectively.
\end{proof}

\subsection{Shapes of subsequences}
\paragraph{Definition of a shape.}
Let $w$ be the word defined in \eqref{def of w}. Suppose that $v\in [k]^m$ is a word for which we wish to bound $M(v,w)$. 

For a map $f\in \mathcal{M}(v,w)$, define  
$$g_i(f)\eqdef\min\{j:f(j)\geq k(i-1)+1\}.$$
We adopt the convention that the minimum is equal to $m+1$ if no $j$ satisfying the condition exists. Note that $k(i-1)+1$ is the position of the first symbol of $\pi_i$ in the definition
of $w$ in \eqref{def of w}.
%(As usual, the minimum of the empty set is $\infty$.)  
%$$g_i(f)\eqdef\begin{cases}\min\{j:f(j)\geq k(i-1)+1\}\\.\end{cases}$$
We write $g_i=g_i(f)$ for brevity. 

Let $h_i$ be the largest integer such that $v\ind{g_i,h_i}$ is a subsequence of $\pi_i$;
such an integer $h_i$ satisfies
\begin{equation}\label{eq:hg}
h_i\geq g_i-1
\end{equation}
because $v\ind{g_i,g_i-1}$ is the empty word.
Set $d_i\eqdef h_i-g_{i+1}$.

To the map $f$ we associate its \emph{shape} $\Shape(f)\in \{0,1,2,3,4,8\}^{n/k-1}$, which is defined
by
$$\Shape(f)_i\eqdef\begin{cases}
    8&\text{if } d_i\geq 10t^4 ,\\
    4&\text{if } 10t^4 >d_i\geq 10t^3 ,\\
    3&\text{if } 10t^3>d_i\geq 10t^2 ,\\
    2&\text{if } 10t^2 >d_i\geq 10t ,\\
    1&\text{if } 10t >d_i\geq 1,\\
    0&\text{if } d_i<1.
\end{cases}$$

\paragraph{The relation between shapes and the number of subsequences.}
We prove an upper bound on $M(w)$ by estimating the number of subsequences of shape $s$, for every possible $s \in \{0,1,2,3,4,8\}^{n/k - 1}$.

As we will see in \Cref{f determined},  the sequence $g_1,g_2,\dotsc,g_{n/k}$ determines the function $f$. To the
first approximation, we approach the task of estimating the number of possible
$g$-sequences sequentially, i.e., given $g_i$, we bound the number of choices for $g_{i+1}$. This is done in \Cref{Simple bound}; it states that the number of choices
is at most~$O(t^{s_i})$. Therefore, we must show that a shape cannot contain too many large elements, which is our next task.

\paragraph{Possible subsequence shapes.} We begin by recording
the basic properties of $g_i$ and $h_i$.
\begin{proposition}The sequences $g_1,g_2,\dotsc$ and $h_1,h_2,\dotsc$ satisfy
\begin{enumerate}[label=\alph*), ref=Proposition \thetheorem(\alph*)]
    \item \label{g monotonicity} $g_i\leq g_{i+1}$ for $i\in [n/k-1],$
   \item \label{gi hi subseq} $v\ind{g_i,h_i}\prec \pi_i$ for $i\in [n/k],$
    \item \label{di nonneg} $d_i\geq -1$ for $i\in [n/k-1],$
    \item \label{f determined} The function $f \in \mathcal{M}(v,w)$ is determined by the sequence $g_1, \dots, g_{n/k}$.
\end{enumerate}
\end{proposition}
\begin{proof}
The properties (a) and (b) are immediate from the definitions of $g_i$ and $h_i$, respectively. 

Suppose, for contradiction's sake, that (c) fails, i.e., $h_i<g_{i+1}-1$. From \eqref{eq:hg}
it follows that $g_i<g_{i+1}$, and so $f(g_i)<ki+1$.
This in turn means that $f(g_i)$ is in the interval $[k(i-1)+1,ki]$, i.e.,
inside the permutation $\pi_i$ in the definition of $w$ in \eqref{def of w}.
The definition of $h_i$ then implies that $f(h_i+1)\geq ki+1$, and so
$h_i\geq g_{i+1}-1$.

 We turn to the proof of (d). If the half-open interval $[g_i,g_{i+1})$ is non-empty, then the preceding argument shows that the $f$-image of the half-open interval $[g_i,g_{i+1})$
lies entirely inside $\pi_i$. Since $\pi_i$ is a permutation, this means that
the values of $f$ on $[g_i,g_{i+1})$ are determined.
\end{proof}

Let $s$ be a shape of some $f\in \mathcal{M}(v,w)$. 
In the next few claims, we show that certain patterns in $s$ cannot occur.
These claims concern the possible values of $s_i,s_{i+1},\dotsc,s_{i+8}$.
To ensure that these are well-defined, we assume that $i<n/k-8$ throughout.

\begin{claim}\label[claim]{b0}
%If $s_i\in \{3,4,8\}$, then either $s_{i+1}=0$ or both $s_{i+1}\in \{1,2\}$ and $s_{i+2}=0.$
  If $s_i\in \{3,4,8\}$, then either 
  \begin{enumerate}[label=\alph*), ref=Proposition \thetheorem(\alph*)]
  \item $s_{i+1}=0$,
  \item $s_{i+1}=1$, and $s_{i+2}=0$,
  \item \label{b0case2} $s_{i+1}=2$, and $s_{i+2}=0$, and $s_{i+j}\in \{0,1\}$ for $j=3,4,\dotsc,7$.
  \end{enumerate}
\end{claim}
\begin{proof}
\textsc{Proof that $s_{i+1}\in \{0,1,2\}$.}
Since the interval $I\eqdef [g_{i+1},\min(h_i,h_{i+1})]$ is a subinterval of both $[g_i,h_i]$ and $[g_{i+1},h_{i+1}]$, it follows that
\( v\ind{I}\prec \pi_i, \pi_{i+1} \).
By \Cref{prop: pi adj} $\LCS(\pi_i, \pi_{i+1}) \leq t^2$, and so $\abs{I}\leq t^2$.
Because $s_3\geq 3$, it follows that $h_i-g_{i+1}=d_i\geq 10t^2$, which means that $I=[g_{i+1},h_{i+1}]$. Thus $h_{i+1}\leq h_i$ and
\[
h_{i+1} - g_{i+1} + 1 \leq t^2.
\]
Therefore,
\[
h_{i+1} - g_{i+2} + 1 \leq t^2.
\]
Hence, \( s_{i+1} \in \{0, 1, 2\} \).\medskip

\textsc{Proof that $s_{i+1}\neq 0$ implies that $s_{i+2}=0$.}
Suppose, towards a contradiction, that
\( h_{i+2} - g_{i+2} \geq 1 \). Because $s_{i+1}\neq 0$, the interval $[g_{i+2},g_{i+2}+1]$ is a subinterval of $[g_{i+1},h_{i+1}]$, and the assumption \( h_{i+2} - g_{i+2} \geq 1 \) implies that  $[g_{i+2},g_{i+2}+1]$ is a subinterval of $[g_{i+2},h_{i+2}]$. Because $h_{i+1}\leq h_i$, the interval $[g_i,h_i]$ contains $[g_{i+1},h_{i+1}]$ and hence also contains $[g_{i+2},g_{i+2}+1]$.
Therefore, \( v\ind{g_{i+2}, g_{i+2}+1}\prec \pi_i,\pi_{i+1}, \pi_{i+2} \). However, this is not possible, for \Cref{prop: pi three adj} asserts that \( \LCS(\pi_i, \pi_{i+1}, \pi_{i+2}) = 1 \). The contradiction shows
that \( h_{i+2} - g_{i+2} \leq 0 \), from which \( s_{i+2} = 0 \) follows.\medskip

\textsc{Proof that $s_{i+2}=2$ implies that $s_{i+j}\in \{0,1\}$ for $j=3,4,\dotsc,7$.}
We claim that
\begin{alignat}{2}\label{b0ind}
  h_{i+1}-g_{i+j}&\geq (12-j)t&&\qquad\text{for }j=2,3,\dotsc,8,\\
    \label{b0mainind}
  h_{i+j}-g_{i+j}&\leq t-1&&\qquad\text{for }j=2,3,\dotsc,7.
%  h_{i+j}-g_{i+j+1}&\leq t&&\qquad\text{for }j=2,3,\dotsc,7.
  \end{alignat}
Note that \eqref{b0mainind} implies that $d_{i+j}\leq h_{i+j}-g_{i+j}\leq t-1$, and thus $s_{i+j}\in \{0,1\}$.
The proof is by induction. The inequality \eqref{b0ind} for $j=2$ follows from $s_{i+1}=2$.

Assuming that \eqref{b0ind} holds for some $j\in \{2,3,\dotsc,7\}$, we will prove that \eqref{b0mainind} holds for the same $j$, and that \eqref{b0ind} holds for $j+1$.
The interval $I_j\eqdef [g_{i+j},\min(h_{i+1},h_{i+j})]$
is contained in both $[g_{i+1},h_{i+1}]=I$ and in $[g_{i+j},h_{i+j}]$, and therefore
$v\ind{g_{i+j},\min(h_{i+1},h_{i+j})}\prec \pi_i,\pi_{i+1},\pi_{i+j}$.
By \Cref{prop: pi three partadj}, \( \LCS(\pi_i, \pi_{i+1}, \pi_{i+j}) \leq t \). In view of \eqref{b0ind},
this implies that $I_j=[g_{i+j},h_{i+j}]$, from which we deduce that \eqref{b0mainind}
holds. 

Therefore,
\begin{alignat*}{2}
  h_{i+1}-g_{i+j+1} &= h_{i+1}-h_{i+j}+d_{i+j} \\
                    &\geq h_{i+1}-h_{i+j}-1&&\qquad\text{by \Cref{di nonneg}}\\
                    &\geq h_{i+1}-g_{i+j}-t,&&\qquad\text{by \eqref{b0mainind}}
\end{alignat*}
proving \eqref{b0ind} for $j+1$.
\end{proof}

\begin{claim}\label[claim]{bb}
    If $s_i=8$, then 
    \begin{enumerate}[label=\alph*), ref=Proposition \thetheorem(\alph*)]
      \item $s_{i+j}\neq 8$ for $j\in [7]$,
      \item\label{bb part b} the interval $[g_{i+j},h_{i+j}]$ is contained in the interval $[g_i,h_i]$ for $j\in [7]$.
    \end{enumerate}
\end{claim}
\begin{proof}
The proof is similar to the proof of the last part of the previous claim.
%We claim, by induction on $j$, that
It is enough to prove that
\begin{alignat}{2}\label{bb 1}
h_i-g_{i+j}\geq (11-j)t^4&&\qquad\text{for }j=1,2,\dotsc,8,\\
\label{bb 2}
h_{i+j}-g_{i+j}\leq t^4-1&&\qquad\text{for }j=1,2,\dotsc,7.
\end{alignat}
Indeed, part (a) of the claim is immediate from \eqref{bb 2} and the monotonicity of the $g$'s, and
part (b) follows from $h_{i+j}\leq h_i$, which can be seen by subtracting \eqref{bb 2} from \eqref{bb 1}.\medskip

The case $j=1$ of the inequality \eqref{bb 1} follows from the assumption $s_i=8$.  Assume that \eqref{bb 1} holds for some value of $j$; we will demonstrate \eqref{bb 2} for the same $j$ as well as \eqref{bb 1} for $j+1$.

The interval $[g_{i+j},\min(h_i,h_{i+j})] $ is a subinterval of both $[g_i,h_i]$ and $[g_{i+j},h_{i+j}]$, and the inequality \( \LCS(\pi_i, \pi_{i+j}) \leq t^4 \) is contained in \Cref{prop: pi nonadj}. In view of \eqref{bb 1},
these two facts imply that $[g_{i+j},\min (h_i,h_{i+j})]=\nobreak[g_{i+j},h_{i+j}]$, from which we deduce \eqref{bb 2}. Combining \Cref{di nonneg} with \eqref{bb 2} yields $h_i-g_{i+j+1}\geq h_i-h_{i+j}-1\geq h_i-g_{i+j}-t^4$, proving \eqref{bb 1} for $j+1$.
%Part (a) is immediate from \eqref{bb 2} and monotonicity of the $g$'s.
%Part (b) follows from $h_{i+j}\leq h_i$, which can be seen by subtracting \eqref{bb 2} from \eqref{bb 1}.
\end{proof}
\begin{claim}\label[claim]{bs0}
    If $s_i=8$ and for some $j\in [6]$, $s_{i+j}\in \{2,3,4\}$, then $s_{i+j+1}\in\{0,1\}$.
\end{claim}
\begin{proof}
   As in the proofs above. Consider an interval $I\eqdef[g_{i+j},\min(h_{i+j},h_{i+j+1})]$. Then $I$ is subinterval of both $[g_{i+j},h_{i+j}]$ and $[g_{i+j+1},h_{i+j+1}]$ and from \Cref{bb part b} it is also a subinterval of $[g_{i},h_{i}]$. Therefore, $v\ind{I}\prec \pi_i,\pi_{i+j},\pi_{i+j+1}$.
   By \Cref{prop: pi three partadj}, $\LCS(\pi_i,\pi_{i+j},\pi_{i+j+1})\leq t$, and hence $|I|\leq t$. Since $h_{i+j}-g_{i+j+1}\geq 10t$, we deduce $I=[g_{i+j+1},h_{i+j+1}]$, from which $s_{i+j+1}\in \{0,1\}$ follows. 
\end{proof}
\begin{claim}\label[claim]{b7}
    If $s_i=8$, then $s_{i+7}\in \{0,1,2\}.$
\end{claim}
\begin{proof}
    From \Cref{bb part b}, $[g_{i+7},h_{i+7}]$ is subinterval of $[g_i,h_i]$, hence $v\ind{g_{i+7},h_{i+7}}\prec \pi_i,\pi_{i+7}$. Therefore, from \Cref{prop: pi adj}, $h_{i+7}-g_{i+7}\leq t^2-1$, from which $d_{i+7}\leq t^2-1$ follows, concluding the proof.
\end{proof}
\begin{claim}\label[claim]{bn}
If $s_i=8$, then among $s_{i+j}$ for $j\in [7]$, there is at most one element in $\{3,4\}$.
\end{claim}
\begin{proof}
    This is similar to the proofs of \Cref{bb} and the last part of \Cref{b0}.

    Let $j$ be the least index such that $s_{i+j}\in \{3,4\}$.
    By \Cref{bb part b}, the interval \( I\eqdef\nobreak [g_{i+j}, h_{i+j}] \) is entirely contained within \( [g_i, h_i] \), and its length $\abs{I}=h_{i+j}-g_{i+j}+1$ satisfies
\[
  h_{i+j}-g_{i+j}\geq h_{i+j} - g_{i+j+1} \geq 10t^2.
\]
With this as the base case, by an inductive argument similar to those in the claims above, applied to the word
$v\ind{g_{i+\ell},\min(h_{i+j},h_{i+\ell})}$, we can show that
\begin{alignat}{2}
h_{i+j}-g_{i+\ell}\geq (11-(l-j))t^2&&\qquad\text{for }\ell=j+1,\dotsc,8,\\
h_{i+\ell}-g_{i+\ell}\leq t^2-1&&\qquad\text{for }\ell=j+1,\dotsc,7,\label{bn final}
\end{alignat}
using the bound $\LCS(\pi_i,\pi_{i+j},\pi_{i+\ell})\leq t^2$ from \Cref{prop:pi distinct}.

As \eqref{bn final} implies that $s_{i+\ell}\notin \{3,4\}$, the claim follows from minimality of~$j$.
\end{proof}

\subsection{Counting \texorpdfstring{$g$-sequences}{g-sequences}}
\paragraph{Basic bound.}
Fix a shape $s$ that satisfies the restrictions enumerated in \Cref{b0} through \Cref{bn} in the preceding
subsection. We shall bound the number of subsequences with shape $s$. Our main tool is
the following bound.

\begin{claim}\label[claim]{Simple bound}
   If $s_i=x$, and $g_i$ is fixed, then there are at most $10t^x$ possibilities for
   the value of $g_{i+1}$.
\end{claim}
\begin{proof}
    If $x\geq 1$, then $g_{i+1}$ is in the interval of length $10t^x-1$, namely $g_{i+1}\in [h_i-10t^{x},h_i-2]$; as $h_i$ is determined by $g_i$, the result follows in this case. If $x=0$, then $g_{i+1}\in \{-1,0\}$.
\end{proof}

\paragraph{Better bound in the case \texorpdfstring{$s_i=8$}{s\textunderscore i=8}.}
Unfortunately, the sequential application of \Cref{Simple bound} by itself is not enough to prove that $M(w)^{k/n}\leq (Ck)^{1/4}$. So, to bridge the gap, we prove that if $s_i=8$,
then some of the few subsequent applications of \Cref{Simple bound} can be strengthened.

For a word $b$ over alphabet $[k]$ and $x\in [8]$, define a set
\[
  E_x(b)\eqdef \{z : b\ind{z}_{\leq x}\neq b\ind{z+1}_{\leq x}\}.
\]

Let $E_{x,y}(b)$ be a subset of $E_x(b)$ that is obtained by starting from the smallest element of $E_x(b)$ and keeping every $y$-th element.
\begin{claim}\label[claim]{b improvement 2}
    If $s_i=8$, and $h_{i+j}-g_{i+j}\geq t^2 y$ for some $j\in [7]$ and some natural number $y$. Then $g_{i+j}$ is at most $t^2y$ away from some element of $E_{5,y}(v)\cap [g_{i+j},h_{i+j}]$.
\end{claim}
\begin{proof}
    Suppose that the interval \( I\eqdef [g_{i+j}, g_{i+j}+t^2y] \) contains no element of \( E_{5,y}(v) \). Then it must contain fewer than \( y \) elements of \( E_5(v) \). 
    Therefore, we can partition \( I \) into at most \( y  \) subintervals, such that every element of $E_5(v)$ is the largest element in its subinterval.

Let $J$ be the longest subinterval in the partition.
    Since $\abs{I}=t^2y+1$, it follows that $\abs{J}\geq t^2+1$.
By \Cref{bb part b}, $v\ind{J}$ is a subsequence of both \( \pi_i \) and \( \pi_{i+j} \). Furthermore, $v\ind{z}_{\leq 5}$ is the same for all $z\in J$.
Hence, from \Cref{prop:prefix 5}, we deduce that
\[
  |J| \leq t^2,
\]
contradicting the choice of~$J$.
\end{proof}
\begin{claim}\label[claim]{b improvement 3}
    If $s_i=8$, and $h_{i+j}-g_{i+j}\geq t^3 y$ for some $j\in [7]$ and some natural number $y$. Then $g_{i+j}$ is at most $t^3y$ away from some element of $E_{3,y}(v)\cap [g_{i+j},h_{i+j}]$.
\end{claim}
\begin{proof}
  The proof is similar to the proof of \Cref{b improvement 2}. However, here we use \Cref{prop:prefix 3} instead. Therefore, the length of the longest subinterval
  in the partition is bounded by~$t^3$.
\end{proof}
\begin{claim}\label[claim]{Monotone}
    Suppose that $b$ is a subsequence of $\pi_i$ for some $i$, and $x\in [8]$. Then $$|E_x(b)|\leq t^x.$$ 
\end{claim}
\begin{proof}
  It is enough to show that $b\ind{y}_{\leq x}$ are distinct for distinct $y\in E_x(b)$.
  Suppose, for contradiction's sake, that $b\ind{y_1}_{\leq x}=b\ind{y_2}_{\leq x}$ for some distinct $y_1,y_2\in E_x(b)$.
  Then the symbols $b\ind{y_1},b\ind{y_1+1},b\ind{y_2}$ cannot appear in this order in any permutation
  of the form $\pi(u)$ because $b\ind{y_1+1}_{\leq x}$ differs from $b\ind{y_1}_{\leq x}=b\ind{y_2}_{\leq x}$.
  In particular, they cannot appear in $\pi_i$, contradicting the assumption on~$b$.
  % TODO: Try being somewhat consistent about using "symbol" rather than "letter" (but "n-symbol word" does not sound right)
\end{proof}

Each permutation of the form $\pi(u)$ (such as $\pi_1,\pi_2,\dotsc$) contains exactly
$t^5$ pairs of adjacent symbols $\ell,\ell'$ such that $\ell_{\leq 5}\neq \ell'_{\leq 5}$.
Therefore, it is reasonable to expect that the gaps between the elements of $E_5(v)$ are on average $t^3$.
Taken together, \Cref{Monotone} and \Cref{b improvement 2} show that, contrary to this expectation,
the gaps between the elements of $E_5(v)$ are at most $t^2$ inside any interval of the form $[g_{i+j},h_{i+j}]$, if $s_i=8$ and $j\in [7]$.
Similarly, \Cref{Monotone} and \Cref{b improvement 3} show a similar result about the gaps in $E_3(y)$.
We use this fact to improve upon \Cref{Simple bound} if $s_i=8$.
\begin{lemma}\label[claim]{Lemm: Bound Big}
If $s_i=8$, and $g_i$ is fixed, and $s_{i+j}\in \{3,4\}$ for some $j\in [6]$, then there are at most $10^5t^{10}$ many possibilities for $(g_{i+1},g_{i+j+1})$.
\end{lemma}
\begin{proof}
  Let \( [g_{i+\ell}, h_{i+\ell}] \) be the longest interval among
  $[g_{i+1},h_{i+1}],[g_{i+2},h_{i+2}],\dotsc,[g_{i+7},h_{i+7}]$.
  Since the length of $[g_{i+j},h_{i+j}]$ exceeds $10t^2$, the
  same is true of the length of $[g_{i+\ell}, h_{i+\ell}]$.

  \cref{bb part b,prop: pi nonadj} together imply that the length of $[g_{i+\ell}, h_{i+\ell}]$ is upper bounded by $t^4$.
  So, there is $q\in [0,\log_2 t-1]$ such that either
\begin{equation}\label{t2 case}
 h_{i+\ell}-g_{i+\ell} \in \left[10 t^2 2^q, 10 t^2 2^{q+1}\right)
\end{equation}
or
\begin{equation}\label{t3 case}
h_{i+\ell}-g_{i+\ell} \in \left[10t^3 2^q, 10 t^3 2^{q+1}\right)
\end{equation}
holds.

Consider the case \eqref{t2 case} first. Assume \( q \) is fixed; set $y\eqdef 5\cdot 2^q$.
By the maximality of \( |[g_{i+\ell}, h_{i+\ell}]| \), and \Cref{di nonneg} it follows that
$g_{i+z+1}-g_{i+z}=h_{i+z}-g_{i+z}-d_{i+z}\leq 4t^2 y$. Therefore,
\[
  g_{i+\ell}-g_{i+z}\leq 24 t^2 y\qquad\text{for all }z\in [7].
\]
Now, from \Cref{b improvement 2}, there exists a point \(p\in E_{5,y}(v)\cap [g_{i+\ell},h_{i+\ell}] \) such that \( g_{i+\ell} \) is within distance \( t^2 y \) of $p$. Hence, both \( g_{i+1} \) and \( g_{i+j+1} \) are within distance at most \( 25t^2 y \) from \( p \).

By \Cref{bb part b}, \( p \in [g_i, h_i] \). Now, $v\ind{g_i,h_i}$ is subsequence of $\pi_i$ and from \Cref{Monotone} we know that
$|E_5(v)\cap [g_i,h_i]|\leq t^5$. Therefore, by the nature of $E_{5,y}$, we deduce that
\[
\abs{E_{5,y}(v)\cap [g_i,h_i]}\leq \frac{t^5}{y}+1.
\]
Hence, for fixed $y$, there are at most $t^5/y+1$ many choices for $p$, and so at most
$(t^5/y+1)\cdot (25t^2 y)^2\leq 2000 t^9y$ many choices for $(g_{i+1},g_{i+j+1})$.

Summing this over all possible $q$ we obtain
\[
  \sum_{q=0}^{\log_2 t-1} 2000 t^9 \cdot (5\cdot 2^q)\leq 10000 t^{10}.
\]

The similar argument in the case \eqref{t3 case} leads to the same bound of $10000t^{10}$ using \Cref{b improvement 3}.
The two cases together imply the stated bound on the number of possibilities for $(g_{i+1},g_{i+j+1})$.
\end{proof}
\begin{lemma}\label[lemma]{clam:701}
    If $s_i=8$, and $g_i$ is fixed, and $s_{i+j}\in\{2,3,4\}$ for some $j\in [3]$, then there are at most $20t^7$ many possibilities for $g_{i+1}$.
\end{lemma}
\begin{proof}
   Let \( \ell \) be the smallest index such that the length of the interval \( [g_{i+\ell}, h_{i+\ell}] \) satisfies
\[
  |[g_{i+\ell}, h_{i+\ell}]| \geq 2t.
\]
Note that $\ell\leq j\leq 3$. By \Cref{prop:prefix 6} the interval \( [g_{i+\ell}, g_{i+\ell}+t+1] \) contains an element of $E_6(v)$, say $z\in E_6(v)\cap [g_{i+\ell}, g_{i+\ell}+t+1]$.

By the same argument as in the preceding lemma, since each interval \( [g_{i+x}, h_{i+x}] \) has length at most \( 2t \) for \( x< j \), it follows that \( g_i \) is within distance at most \( 10t \) from \( z \). Since $v\ind{g_i,h_i}\prec \pi_i$, from \Cref{Monotone} it follows that there are $t^6$ many choices for~$z$,
and so the number of possibilities for $g_{i+1}$ is at most $20t^7$.
\end{proof}
\begin{lemma}\label[lemma]{lemma 1}
  If $s_i=8$, and $g_i$ is fixed. 
  \begin{enumerate}[label=\alph*), ref=Lemma \thetheorem(\alph*)]
    \item If $(s_{i+6},s_{i+7})\in\{(3,1),(4,1)\}$, then there are at most $10^{18}t^{17}$ many possibilities for the $9$-tuple $(g_{i+1},\dots g_{i+9})$.
    \item Otherwise, there are at most $(10t)^{16}$ many possibilities for the $8$-tuple $(g_{i+1},\dots g_{i+8})$.
  \end{enumerate}
\end{lemma}
\begin{proof}
   \textbf{Case 1:} \( s_{i+1} = 2 \)

   From \Cref{b0}, we know that \( s_{i+j} \in \{0, 1\} \) for all \( j \in [2, 7] \). In particular, since $s_{i+6}\in \{0,1\}$, we are in case (b) of the Lemma.
   By \Cref{Simple bound}, given $g_{i+j}$, the number of possibilities for $g_{i+j+1}$ is at most $10t^{s_i}$. Using this bound for $j=0,1,\dotsc,7$
   iteratively, we see that the number of possibilities for $(g_{i+1},\dotsc,g_{i+8})$, for fixed $g_i$, is at most $10^8t^{\sum_{j=0}^7 s_i}\leq 10^8t^{16}$.

\vspace{0.5em}
\textbf{Case 2:} There exists \( j \in [7] \) such that \( s_{i+j} \in \{3, 4\} \), and $(s_{i+6},s_{i+7})\notin\{(3,1),(4,1)\}$.

%By \Cref{b7}, $j\neq 7$. By this lemma's assumption, $j\neq 6$. So, $j\in [5]$.

From \Cref{Lemm: Bound Big}, there are at most \( 10^5  t^{10} \) choices for the pair \( (g_{i+1},g_{i+j+1}) \).
Fixing $(g_{i+1},g_{i+j+1})$, we bound the number of possibilities for \( (g_{i+2}, \dots, \widehat{g_{i+j+1}}, \dots, g_{i+8}) \).
By \Cref{Simple bound}, given $g_{i+\ell}$, the number of possibilities for $g_{i+\ell+1}$ is at most $10t^{s_{i+\ell}}$. Using this bound for $\ell=1,2,\dotsc,\hat{j},\dotsc,7$,
we see that the number of possibilities for \( (g_{i+2}, \dots, \widehat{g_{i+j+1}}, \dots, g_{i+8}) \) is at most
\[
  10^6\prod_{\substack{\ell\in [7]\\\ell\neq j}} t^{s_{i+\ell}}.
\]
We thus need to bound the sum of the sequence $s_{i+1},s_{i+2},\dotsc,\widehat{s_{i+j}},\dotsc,s_{i+7}$. For an integer $r$, let $C_r$ be the number
of elements in this sequence that are at most~$r$.
\begin{enumerate}[label=\Roman*)]
\item From \Cref{b0} applied at positions $i$ and $i+j$, we know that both $\{s_{i+1},s_{i+2}\}$ and $\{s_{i+j+1},s_{i+j+2}\}$ contain $0$.
  If $j\leq 5$, this implies that $C_0\geq 2$.
  Note $j\neq 7$ by \Cref{b7}.
  If $j=6$, then the assumption $(s_{i+6},s_{i+7})\notin\{(3,1),(4,1)\}$ together with \Cref{b0} and \Cref{bs0} imply that $s_7=0$, and so $C_0\geq 2$
  in this case as well.  
\item Since $s_i=8$ and $s_{i+j}\in \{3,4\}$, \Cref{b0} implies that $s_{i+1}\in \{0,1\}$. By \Cref{bs0} applied at positions $i+2,i+4$ and $i+6$,
we know that at least three elements of $s_{i+2},\dotsc,s_{i+7}$ are in $\{0,1\}$, and so $C_1\geq 4$.
\item From \Cref{bn}, we know that $C_2=6$.
\end{enumerate}
Therefore
\[
  \sum_{\substack{\ell\in [7]\\\ell\neq j}} s_{i+\ell}=(6-C_0)+(6-C_1)+(6-C_2)\leq 6,
\]
and so the number of possibilities for $(g_{i+1},\dotsc,g_{i+8})$ is at most $10^5 t^{10}\cdot 10^6t^6=10^{11}t^{16}$.

\vspace{0.5em}
\textbf{Case 3:} We have $(s_{i+6},s_{i+7})\in\{(3,1),(4,1)\}$.

Similarly to the previous case, we bound the number of possibilities for $(g_{i+1},g_{i+7})$ by $10^5t^{10}$,
and the number of possibilities for \( (g_{i+2}, \dots, \widehat{g_{i+j+1}}, \dots, g_{i+9}) \) by
\[
  10^7\prod_{\substack{\ell\in [8]\\\ell\neq j}} t^{s_{i+\ell}}.
\]
Since the assumption $(s_{i+6},s_{i+7})\in\{(3,1),(4,1)\}$ and \Cref{b0}
imply that $s_{i+8}=0$, the proof reduces to the task of bounding the sum of $s_{i+1},s_{i+2},\dotsc,\widehat{s_{i+j}},\dotsc,s_{i+7}$.
This is done the same way as in Case~2.

\vspace{0.5em}
\textbf{Case 4:} There exists \( j \in [3] \) such that \( s_{i+j} = 2 \), and \( s_{i+\ell} \not\in \{3,4\} \) for all \( \ell \in [7] \), and \( s_{i+1} \neq 2 \).

Since none of the $s_{i+\ell}$ exceed $2$, by \Cref{bs0} implies that $s_{i+\ell}+s_{i+\ell+1}\leq 3$ for any $\ell\in [6]$.
Hence,
$s_{i+4}+s_{i+5}+s_{i+6}+s_{i+7}\leq  6$.
% $\sum_{\ell=4}^7 s_{i+\ell}\leq  6$.
This also implies that, if $s_{i+1}=0$, then
\begin{equation}\label{eq:ell three}
  s_{i+1}+s_{i+2}+s_{i+3}\leq 3.
%  \sum_{\ell=1}^3 s_{i+\ell}\leq 3.
\end{equation}
On the other hand, if $s_{i+1}=1$, then \eqref{eq:ell three} holds by \Cref{b0}. So, in either case
%$\sum_{\ell=1}^7 s_{i+\ell}\leq 9$.
$s_{i+1}+\dotsb+s_{i+7}\leq 9$.
By \Cref{clam:701}, there are at most \( 20 t^7 \) possibilities for \( g_{i+1} \), for fixed $g_i$.
So, the total number of possibilities for $(g_{i+1},\dotsc,g_{i+8})$ is at most $20t^7\cdot 10^7t^9\leq 10^9t^{16}$.

\vspace{0.5em}
\textbf{Case 5:}
% $s_{i+1},s_{i+2},s_{i+3}\neq 2$, and \( s_{i+\ell} \not\in \{3,4\} \) for \( \ell \in [7] \).
We have $s_{i+j}\neq 2$ for \( j \in [3] \), and \( s_{i+\ell} \not\in \{3,4\} \) for \( \ell \in [7] \).

Since none of $s_{i+1},s_{i+2},s_{i+3}$ are equal to $2$, by \Cref{b0}
\begin{equation}\label{eq:ell two}
  s_{i+1}+s_{i+2}+s_{i+3}\leq 2.
\end{equation}
The rest of the proof is then the same as in Case 3, using \Cref{Simple bound} and \eqref{eq:ell two} instead of \Cref{clam:701} and \eqref{eq:ell three}, respectively.
\end{proof}

\subsection{Putting everything together (proof of the upper bound in \texorpdfstring{\Cref{Very Main}}{Theorem 2})}
Define \( S \eqdef \{0,1,2,3,4,8\}^{n/k-1} \). Then 
\[
M(v, w) = \sum_{s \in S} |\Shape^{-1}(s)|.
\]

Fix \( s \in S \). We decompose \( s \) into segments inductively as follows.
Suppose that the prefix $s_1s_2\dots s_{i-1}$ has been partitioned into the segments already.
Depending on the value of $(s_i,s_{i+1},\dotsc,s_{i+7})$, we add a new segment according to the following rule:
\begin{itemize}
\item If \( s_i \in \{0,1,2\} \), let the singleton \( s_i\) be the next segment.
\item If \( s_i \in \{3,4\} \) and \( s_{i+1} = 0 \), let \( s_i s_{i+1} \) be the next segment.
\item If \( s_i \in \{3,4\} \) and \( s_{i+1} \in \{1,2\} \), let \( s_i s_{i+1} s_{i+2} \) be the next segment.
%\item If \( s_i = 8 \), then if \( s_{i+6} \in \{3,4\} \) and \( s_{i+7} = 1 \), let \( s_i \dots s_{i+8} \) be the next segment. Otherwise, let \( s_i \dots s_{i+7} \) be the next segment.
\item If \( s_i = 8 \) and \( (s_{i+6},s_{i+7}) \in \{(3,1),(4,1)\} \), let \( s_i \dots s_{i+8} \) be the next segment.
\item If \( s_i = 8 \) and \( (s_{i+6},s_{i+7}) \notin \{(3,1),(4,1)\} \), let \( s_i \dots s_{i+7} \) be the next segment.
\end{itemize}
We repeat the process until fewer than $9$ elements are left.

From \Cref{lemma 1} and \Cref{b0} together with \Cref{Simple bound}, we conclude that for each just-defined segment of length \( x \), say $s_is_{i+1}\dots s_{i+x-1}$ we have at most \( (100 t^2)^x \) many choices for the corresponding values of \( g_{i+1},g_{i+2},\dots g_{i+x} \), once $g_{i}$ is fixed.

Since these segments account for all of $s$, except possibly at most $8$ terminal symbols, it follows from \Cref{Simple bound} that
\[
|\Shape^{-1}(s)| \leq (100 t^2 )^{n/k} \cdot (10t^8)^8.
\]

Since \( |S| = 6^{n/k - 1} \) and $n\geq Ck\log k$, for large enough $C$, we obtain
\[
M(v, w) \leq 6^{n/k} \cdot (100 t^2)^{n/k}\cdot (10t^8)^8 \leq(C k)^{n/(4k)}.
\]
\bibliographystyle{plain}
\bibliography{refs.bib}

\end{document}